\documentclass[preprint,12pt]{elsarticle}

\usepackage{amsmath,amssymb}
\newcommand{\CC}{\mathbb{C}}

\newcommand{\EE}{\mathcal{E}}
\newcommand{\rme}{\mathrm{e}}
\newcommand{\rmi}{\mathrm{i}}
\newcommand{\rmd}{\mathrm{d}}
\newcommand{\Fra}[2]{\displaystyle \frac{ #1}{ #2}}
\newcommand{\Lapl}[1]{\displaystyle {\nabla}^2 #1}
\newcommand{\abs}[1]{\left\lvert#1\right\rvert}

\newcommand{\de}{\,\mathrm{d}}

\newcommand{\RR}{\mathbb R}
\newcommand{\matlab}{Matlab$^\text{\textregistered}$ R2014b}
\newcommand{\octave}{GNU Octave 4.0.0}
\newcounter{bla}

\begin{document}

\begin{frontmatter}



\title{INFFTM: Fast evaluation of 3d Fourier series
in MATLAB with an application to\\quantum vortex reconnections}


\author[a]{Marco Caliari\corref{author}}
\author[a]{Simone Zuccher}

\cortext[author] {Corresponding author.\\\textit{E-mail address:} marco.caliari@univr.it}
\address[a]{Department of Computer Science, University of Verona (Italy)}

\begin{abstract}
Although Fourier series approximation is ubiquitous in computational physics 
owing to the Fast Fourier Transform (FFT) algorithm,
efficient techniques for the fast evaluation of a 
three-dimensional truncated Fourier series at a set of \emph{arbitrary} 
points are quite rare, especially in MATLAB language.
Here we employ the Nonequispaced Fast Fourier Transform (NFFT, by 
J.~Keiner, S.~Kunis, and D.~Potts), 
a C~library designed for this purpose, and provide a 
Matlab$^\text{\textregistered}$ and GNU Octave
interface that makes NFFT easily available to the Numerical Analysis community.
We test the effectiveness of our package in the framework of quantum vortex
reconnections, where pseudospectral Fourier methods are commonly used and 
\emph{local} high resolution is required in the post-processing stage.
We show that the efficient evaluation of a truncated Fourier series at arbitrary
points provides excellent results %can be performed efficiently and 
at a computational cost much smaller than 
carrying out a numerical simulation of the problem on a sufficiently fine
regular grid that can reproduce comparable details of the reconnecting 
vortices.

%In addition to the manuscript you must supply: the program source code; job control scripts, where applicable; a README file giving the names and a brief description of all the files that make up the package and clear instructions on the installation and execution of the program; sample input and output data for at least one comprehensive test run; and, where appropriate, a user manual. These should be sent, via email as a compressed archive file, to the CPC Program Librarian at cpc@qub.ac.uk.
\end{abstract}


\end{frontmatter}

% \linenumbers

% Computer program descriptions should contain the following
% PROGRAM SUMMARY.

\section{Introduction}
Fourier series approximation is a fundamental tool in computational physics.
The main reason for its widespread usage is the availability of the Fast
Fourier Transform (FFT) algorithm which allows
to evaluate, in the three-dimensional case, 
a linear combination of $N_1N_2N_3$ trigonometric
polynomials at a sample of $N_1N_2N_3$ points
of a regular grid
with a computational cost $\mathcal{O}(N_1N_2N_3(\log N_1+\log N_2+\log N_3))$ 
instead
of the cost $\mathcal{O}(N_1^2N_2^2N_3^2)$ 
of a direct Discrete Fourier Transform.

FFT, as implemented in the FFTW \cite{FFTW05} library, is nowadays available 
and easy to use
on most high level computational tools. For instance, in the 
MATLAB\footnote{We will refer to MATLAB
as the programming language used by the softwares 
Matlab$^\text{\textregistered}$ and 
GNU Octave.}
language the functions
\verb+ifft+, \verb+ifft2+, and \verb+ifftn+ allow fast evaluation of
trigonometric polynomials at a specific regular grid of points in one-, two-
and $n$-dimensions, respectively. 
The fast evaluation of a three-dimensional truncated Fourier series
at a set of \emph{arbitrary} points is a more challenging task.
NFFT \cite{KKP09} 
is a C library that approximates the evaluation of 
a truncated Fourier series at a set
of $M$ arbitrary points at cost 
$\mathcal{O}(N_1N_2N_3(\log N_1+\log N_2+\log N_3)+
M\lvert \log \varepsilon\rvert^3)$, where $\varepsilon$ is the desired
accuracy (typically the double precision one).
Although a Matlab$^\text{\textregistered}$ interface is provided,
to our knowledge it is not commonly used in computational science, maybe 
because of the NFFT algorithm formulation.
%to the one-dimensional case for sake of clarity, 
For instance, given the set of coefficients $\{\hat\psi_k\}_{k=1}^N$, 
\verb+ifft+ performs the fast evaluation of
\begin{equation*}
\frac{1}{N}\sum_{k=1}^N\hat \psi_k \rme^{2\pi\rmi(k-1)x_n},
\quad x_n=\frac{n-1}{N},
\quad n=1,2,\ldots,N ,
\end{equation*}
whereas, given the set of coefficients $\{\hat\psi_k\}_{k=-N/2}^{N/2-1}$, 
NFFT evaluates 
\begin{equation*}
\sum_{k=-N/2}^{N/2-1}\hat \psi_k\rme^{-2\pi\rmi k \zeta_m},
\quad \zeta_m\in\left[-\frac{1}{2},\frac{1}{2}\right),\quad m=1,2,\ldots,M
\end{equation*}
in an approximated and fast way.

The present work provides the package INFFTM, a MATLAB interface based on NFFT
for the fast evaluation of a truncated Fourier series of a function 
$\psi\colon\Omega\to\CC$ at a set $\Xi_M$ of arbitrary $M$ points in 
the computational domain $\Omega=\prod_{d=1}^3[a_d,b_d)$.
%The package is named INFFTM and its main function is \verb+infft3+, which is 
%invoked by the simple command
%\begin{verbatim}
%infft3(psihat,a,b,Xi)
%\end{verbatim}
The intermediate case of the evaluation of a truncated Fourier series 
at an arbitrary rectilinear grid $Y_{\boldsymbol M}$ with $M_1M_2M_3$ points 
is also addressed. 
%The function \verb+igridftn+ had been specifically designed for this purpose and
%works in any space dimension. 
%The evaluation of truncated Fourier series at rectilinear
%grids of arbitrary points is fully described in~\S~\ref{sec:fourier}.

The motivation for developing this tool relies on the need for \emph{localized} 
high resolution encountered in the post-processing stage of reconnecting 
quantum vortices~\cite{ZCBB12,ZR15}.
The dynamics of quantum vortices and their possible reconnections are
properly described by the Gross--Pitaevskii equation, which is normally solved
by resorting to the Time Splitting pseudoSPectral (TSSP) approach.
%, which
%solves one part of the equation in the Fourier space and the other part
%in the physical space. 
Since the details of reconnections are localized in space at scales 
much smaller than the vortex core size, employing standard FFT on
regular grids would require an excessive, and thus infeasible, number of 
grid points in order to achieve the required resolution. 
However, in~\cite{CZ16x} it was found that, even in the presence of singular 
solutions, the number of Fourier coefficients required for an accurate 
description of the solution is not very large.
The need for high local resolution at the post-processing stage, however, 
urgently demanded a tool for the efficient evaluation of such 
a \emph{small} truncated Fourier series
at a localized set of clustered points.

The paper is organized as follows.
In~\S~\ref{sec:fourier} we present the details of Fourier series decomposition 
and evaluation in three dimensions, whereas in~\S~\ref{sec:qvr} we describe 
the framework for quantum fluids simulations.
In~\S~\ref{sec:prog} we outline the main functions of our INFFTM package 
and in~\S~\ref{sec:numexp} we show the result of the two main drivers 
performing a quantum vortex reconnection and some evaluations of the truncated
Fourier series at different rectilinear grids and arbitrary points.
\section{Fourier series decomposition and evaluation}\label{sec:fourier}
This section, which is the core of the whole work, introduces the necessary
notation and describes how the Fourier decomposition is performed, together
its successive evaluation at rectilinear grids or arbitrary points.
Let be
\begin{equation*}
I_{\boldsymbol N}=\prod_{d=1}^3\{1,2,\ldots,N_d\},\quad \boldsymbol N=(N_1,N_2,N_3),\ 
\text{$N_d$ even}.
\end{equation*}
Given a complex function $\psi\in L^2(\Omega)$, with 
$\Omega=\prod_{d=1}^3[a_d,b_d)$, 
its truncated Fourier series is
\begin{equation}\label{eq:fs}
\hat \psi(\boldsymbol x)=
\sum_{\boldsymbol k\in I_{\boldsymbol N}}\hat \psi_{\boldsymbol k}\EE_{\boldsymbol k}(\boldsymbol x),
\quad \hat\psi_{\boldsymbol k}\in\CC,
\end{equation}
where $\boldsymbol x=(x_1,x_2,x_3)\in\Omega$, $\boldsymbol k=(k_1,k_2,k_3)$ 
is a multiindex and 
\begin{equation*}
\EE_{\boldsymbol k}(\boldsymbol x)=
\prod_{d=1}^3\frac{\rme^{2\pi\rmi (k_d-1-N_d/2)(x_d-a_d)/(b_d-a_d)}}
{\sqrt{b_d-a_d}} .
\end{equation*}
Given the regular grid of points
\begin{equation*}
X_{\boldsymbol N}=\{\boldsymbol x_{\boldsymbol n}=(x_{1,n_1},x_{2,n_2},x_{3,n_3})\}=
\prod_{d=1}^3\{a_d+(n_d-1)h_d,\ n_d=1,2,\ldots,N_d\} ,
\end{equation*}
with $h_d=(b_d-a_d)/N_d$,
the approximate Fourier 
coefficients $\hat \psi_{\boldsymbol k}$ are computed by the 
three-dimensional trapezoidal quadrature formula applied to the integral
\begin{equation*}
\int_{\Omega} \psi(\boldsymbol x)\overline{\mathcal E}_{\boldsymbol k}(\boldsymbol x)
\rmd \boldsymbol x ,
\end{equation*}
where $\overline{\mathcal E}_{\boldsymbol k}(\boldsymbol x)$ denotes the 
complex conjugate
of $\mathcal E_{\boldsymbol k}(\boldsymbol x)$. The function 
$\hat \psi(\boldsymbol x)$ turns out to be an approximation of the original
$\psi(\boldsymbol x)$ which interpolates it at the points $X_{\boldsymbol N}$.
The denominator in the basis functions $\EE_{\boldsymbol k}$ assures
the equivalence
\begin{equation*}
\sum_{\boldsymbol k\in I_{\boldsymbol N}}\lvert \hat\psi_{\boldsymbol k}\rvert^2=
\int_{\Omega}\lvert \hat\psi(\boldsymbol x)\rvert^2\rmd \boldsymbol x=
h_1h_2h_3\sum_{\boldsymbol x_{\boldsymbol n}\in X_{\boldsymbol N}}
\lvert\hat \psi(\boldsymbol x_{\boldsymbol n})\rvert^2\approx
\int_{\Omega}\lvert \psi(\boldsymbol x)\rvert^2\rmd \boldsymbol x
\end{equation*}
for any domain $\Omega=\prod_{d=1}^3[a_d,b_d)$.

The regular grid of points $X_{\boldsymbol N}$
can be represented in MATLAB by
\begin{verbatim}
[X{1:3}] = ndgrid(x{1:3})
\end{verbatim}
where
\begin{verbatim}
x{d} = linspace(a(d),b(d),N(d)+1)'; x{d} = x{d}(1:N(d))
\end{verbatim}
Given $\boldsymbol n=\verb+[n(1),n(2),n(3)]+$, we have
\begin{equation*}
\begin{aligned}
\boldsymbol x_{\boldsymbol n}&=\verb+[x{1}(n(1)),x{2}(n(2)),x{3}(n(3))]+=\\
&=\verb+[X{1}(n(1),n(2),n(3)),...+\\
&\phantom{{}={}}\verb+ X{2}(n(1),n(2),n(3)),...+\\
&\phantom{{}={}}\verb+ X{3}(n(1),n(2),n(3))]+
\end{aligned}
\end{equation*}
If \verb+psi+ denotes the MATLAB three-dimensional array containing the
values of $\psi$ at $X_{\boldsymbol N}$, then the three-dimensional array
\verb+psihat+ of approximate Fourier coefficients $\hat \psi_{\boldsymbol k}$ is 
recovered using the fast Fourier 
transform
\begin{verbatim}
psihat = fftshift(fftn(psi)) * prod(sqrt(b - a) ./ N)
\end{verbatim}
whose computational cost is $\mathcal{O}(N_1N_2N_3(\log N_1+\log N_2+\log
N_3))$. 

%If the grid is given in terms of \verb+meshgrid(x{1:3})+, then the
%corresponding values of $\psi$ are \verb+permute(psi,[2,1,3])+.

Given the truncated 
Fourier series approximation of a function, it is trivial to 
approximate
its partial derivatives with respect to the directions $x_d$ since
%as the partial derivative of a basis function with respect to $x_d$ is
\begin{equation}\label{eq:Lambda}
\partial_{x_d} \EE_{\boldsymbol k}(\boldsymbol x)=\Lambda_{\boldsymbol k,d}\EE_{\boldsymbol k}(\boldsymbol x),\quad
\Lambda_{\boldsymbol k,d}=2\pi\rmi(k_d-1-N_d/2)/(b_d-a_d),
\end{equation}
which leads to
\begin{equation*}
\partial_{x_d}\hat \psi(\boldsymbol x)=
\sum_{\boldsymbol k\in I_{\boldsymbol N}}\Lambda_{\boldsymbol k,d}\hat \psi_{\boldsymbol k}
\EE_{\boldsymbol k}(\boldsymbol x) .
\end{equation*}
From the definition of $\Lambda_{\boldsymbol k,d}$, it follows that
\begin{equation*}
\begin{split}
\int_\Omega\lvert\nabla \hat \psi(\boldsymbol x)\rvert^2\rmd \boldsymbol x&=
\sum_{\boldsymbol k\in I_{\boldsymbol N}}\left(\abs{\Lambda_{\boldsymbol k,1}}^2+
\abs{\Lambda_{\boldsymbol k,2}}^2+
\abs{\Lambda_{\boldsymbol k,3}}^2\right)\lvert\hat \psi_{\boldsymbol k}\rvert^2=\\
&=-\sum_{\boldsymbol k\in I_{\boldsymbol N}}\left(\Lambda_{\boldsymbol k,1}^2+
\Lambda_{\boldsymbol k,2}^2+
\Lambda_{\boldsymbol k,3}^2\right)\lvert\hat \psi_{\boldsymbol k}\rvert^2
=-\int_\Omega \nabla^2\hat \psi(\boldsymbol x)\overline{\hat\psi(\boldsymbol x)}
\rmd \boldsymbol x,
\end{split}
\end{equation*}
where the last equivalence comes from integration by parts and taking into
account the periodicity of $\hat \psi(\boldsymbol x)$ in the 
computational domain $\Omega$.
\subsection{Evaluation of a truncated Fourier series at a rectilinear grid}
The evaluation of a truncated 
Fourier series at the regular grid $X_{\boldsymbol N}$ 
can be implemented straightforwardly by employing the inverse fast Fourier 
transform
\begin{verbatim}
psihathat = ifftn(ifftshift(psihat)) / prod(sqrt(b - a) ./ N)
\end{verbatim}
whose computational cost is 
$\mathcal{O}(N_1N_2N_3(\log N_1+\log N_2+\log N_3))$.

Given an arbitrary rectilinear grid $Y_{\boldsymbol M}=\prod_{d=1}^3\{y_{d,m_d},\ m_d=1,2,\ldots, M_d\}\subset\Omega$,
we introduce the matrices
\begin{equation*}
\EE^d=\left(e^d_{m_dk_d}\right)=\frac{\rme^{2\pi\rmi (k_d-1-N_d/2)(y_{d,n_d}-a_d)/(b_d-a_d)}}{\sqrt{b_d-a_d}}\in\CC^{M_d\times N_d},
\quad d=1,2,3
\end{equation*}
and then, for $\boldsymbol y_{\boldsymbol m}\in Y_{\boldsymbol M}$, evaluate
\begin{equation*}
\hat\psi(\boldsymbol y_{\boldsymbol m})=\sum_{\boldsymbol k\in I_{\boldsymbol N}}
\hat\psi_{\boldsymbol k}\EE_{\boldsymbol k}(\boldsymbol y_{\boldsymbol m})=
\sum_{k_3=1}^{N_3}e^3_{m_3k_3}\left(\sum_{k_1=1}^{N_1}e^1_{m_1k_1}\left(
\sum_{k_2=1}^{N_2}\hat \psi_{(k_1,k_2,k_3)}e^2_{m_2 k_2}\right)\right).
\end{equation*}
We observe that the inner sum corresponds, for each $k_3$, to a matrix-matrix
product between $\hat \psi_{(k_1,k_2,k_3)}$ and the matrix $\EE^2$ transposed,
leading to a computational cost $\mathcal{O}(N_1N_2M_2)$. 
The middle sum corresponds to a second matrix-matrix product between the 
previous result and the matrix $\EE^1$, for computational cost 
$\mathcal{O}(N_1M_2M_3)$. 
If the result, i.e. a matrix of order $M_1\times M_2$, is computed and 
stored for each $k_3$ 
(this cost is $\mathcal{O}(N_3(N_1N_2M_2+N_1M_2M_3))$), then the outer sum 
corresponds to the multiplication of the
term $e^3_{m_3k_3}$ for such matrices, for a total computational cost 
$\mathcal{O}(N_3M_1M_2M_3)$.
A straightforward implementation in MATLAB of this strategy could be
\begin{verbatim}
for d = 1:3
  E{d} = exp(2*pi*1i * (y{d} - a(d)) / (b(d) - a(d)) ...
             * (-N(d)/2:N(d)/2 - 1)) / sqrt(b(d) - a(d));
end
psihaty = zeros(M);
for k3 = 1:N(3)
  temp = E{1} * (psihat(:,:,k3) * E{2}.'); % 2d evaluation   
  for m3 = 1:M(3)
    psihaty(:,:,m3) = psihaty(:,:,m3) + temp * E{3}(m3,k3);
  end
end
\end{verbatim}
The routine \verb+ndcovlt+\footnote{It was originally written by 
Jaroslav Hajek for the linear-algebra package of GNU Octave.} implements the 
same evaluation avoiding
the two loops over $N_3$ and $M_3$. Both implementations do not need the 
explicit construction of $Y_{\boldsymbol M}$. 
In order to have an idea of the computational cost, we tested the evaluation 
of a series with $64^3$ random complex Fourier coefficients
at the regular grid $X_{\boldsymbol M}$ with $64^3$ points and obtained what
follows.
\begin{verbatim}
ifft
Elapsed time is 0.010371 seconds.
ndcovlt

error_inf =

   1.0725e-13

Elapsed time is 0.049872 seconds.
two loops

error_inf =

   1.0760e-13

Elapsed time is 0.201674 seconds.
\end{verbatim}
in \matlab{} and
\begin{verbatim}
ifft
Elapsed time is 0.0108922 seconds.
ndcovlt
error_inf =    9.5313e-14
Elapsed time is 0.0321529 seconds.
two loops
error_inf =    9.4936e-14
Elapsed time is 0.258168 seconds.
\end{verbatim}
in \octave{}. The first elapsed time is due to the inverse fast Fourier
transform whose result is used to measure the error, in infinity norm,
with respect to the other two methods. In the other two methods,
the computational cost for the evaluation of $\mathcal{E}^d$, $d=1,2,3$,
is not considered. This test can be found at the end of the \verb+igridftn.m+ 
file
and can be run, in GNU Octave, by \verb+demo igridftn+. 
Due to the randomness of the Fourier coefficients the values of the errors
are not perfectly reproducible.
The implementation via \verb+ndcovlt+ is always much faster than the usage of
nested loops; the factor is four in \matlab{} and eight in \octave{}, 
where JIT (Just-in-time accelerator) is not available. 
Be observe that a truncated 
Fourier series can be evaluated at any rectilinear grid 
$Y_{\boldsymbol M}$ and that \verb+ndcovlt+ is very general as it can evaluate 
an $n$-dimensional truncated series at a rectilinear grid of points. 
For instance, it was used in~\cite{CR13} for the evaluation of 
truncated Hermite series.
\subsection{Evaluation of a truncated Fourier series at arbitrary points}
Given a set of arbitrary points $\Xi_M=\{\boldsymbol \xi_m=({\xi_1}_m,{\xi_2}_m,{\xi_3}_m),\ m=1,2,\ldots, M\}\subset\Omega$, it is possible to evaluate 
$\hat\psi(\boldsymbol \xi_m)$ by firstly computing
\begin{equation*}
\EE^d=\left(e^d_{k_d}\right)\frac{\rme^{2\pi\rmi(k_d-1-N_d/2)({\xi_d}_m-a_d)/(b_d-a_d)}}{\sqrt{b_d-a_d}}
\in\CC^{N_d},\quad d=1,2,3
\end{equation*}
and then
\begin{equation}\label{eq:Fs}
\hat \psi(\boldsymbol \xi_m)=\sum_{\boldsymbol k\in I_{\boldsymbol N}}
\hat\psi_{\boldsymbol k}\EE_{\boldsymbol k}(\boldsymbol \xi_m)=
\sum_{k_1=1}^{N_1}\sum_{k_2=1}^{N_2}\sum_{k_3=1}^{N_3}\hat\psi_{(k_1,k_2,k_3)}e^1_{k_1}
e^2_{k_2}e^3_{k_3}.
\end{equation}
This can be done in MATLAB by the code
\begin{verbatim}
for m = 1:M
  E{1}(:,1,1) = exp(2*pi*1i * (Xi(1,m) - a(1)) / (b(1) - a(1)) ...
                    * (-N(1)/2:N(1)/2 - 1)) / sqrt(b(1) - a(1));
  E{2}(1,:,1) = exp(2*pi*1i * (Xi(2,m) - a(2)) / (b(2) - a(2)) ...
                    * (-N(2)/2:N(2)/2 - 1)) / sqrt(b(2) - a(2));
  E{3}(1,1,:) = exp(2*pi*1i * (Xi(3,m) - a(3)) / (b(3) - a(3)) ...
                    * (-N(3)/2:N(3)/2 - 1)) / sqrt(b(3) - a(3));
  EE = bsxfun(@times,E{1} * E{2},E{3});
  psihatxi(m) = sum(psihat(:) .* EE(:));
end
\end{verbatim}
at a computational cost $\mathcal{O}(MN_1N_2N_3)$. 
This implementation is limited to three dimensions, but it can be extended 
to any $n$-dimensional truncated series. Unfortunately, due to the 
construction of the vectors $\mathcal{E}^d$ inside a loop, this implementation
turns out to be quite inefficient.
\subsubsection{NFFT}
Provided the set of points
$\{\boldsymbol \zeta_m=({\zeta_1}_m,{\zeta_2}_m,{\zeta_3}_m),\
m=1,2,\ldots,M\}$, with $-1/2\le {\zeta_d}_m< 1/2$, 
NFFT performs a fast approximation of
\begin{equation*}
\hat f(\boldsymbol \zeta_m)=\sum_{\boldsymbol k\in I_{\boldsymbol N}}\left(
\hat f_{\boldsymbol k}
\prod_{d=1}^3 \rme^{-2\pi\rmi (k_d-1-N_d/2){\zeta_d}_m}\right).
\end{equation*}
%\begin{equation*}
%\hat f(\zeta_m)=\sum_{k=1}^{N}\hat f_{k}\rme^{-2\pi\rmi (k-1-N/2){\zeta}_m},\quad 1\le m\le M,\quad -\frac{1}{2}\le {\zeta}_m<\frac{1}{2}
%\end{equation*}
%If we define
%\begin{equation*}
%\zeta_m=-\left(\frac{\xi_m-a}{b-a}-\frac{1}{2}\right)
%\end{equation*}
%and
%\begin{equation*}
%\hat f_{k}=\hat\psi_{k}\frac{\rme^{\pi\rmi (k-1-N/2)}}{\sqrt{b-a}}
%\end{equation*}
%we have
%\begin{equation*}
%\hat f(\zeta_m)=\sum_{k=1}^{N}\hat \psi_{k}\frac{\rme^{\pi\rmi (k-1-N/2)}\rme^{2\pi\rmi (k-1-N/2)(\xi_m-a)/(b-a)}\rme^{-\pi\rmi (k-1-N/2)}}{\sqrt{b-a}}=
%\hat\psi(\xi_m)
%\end{equation*}
%We notice that $\rme^{\pi\rmi (k-1-N/2)}=(-1)^{k-1-N/2}$, $k=1,2,\ldots,N$. 
%The extension to the three-dimensional case is not difficult.
Given the coefficients $\{\hat \psi_{\boldsymbol k}\}_{\boldsymbol k}$ 
and the evaluation points $\{\boldsymbol \xi_m\}_m$, Fourier 
series~\eqref{eq:fs} evaluation at $\Xi_M$ 
can be approximated by calling the NFFT algorithm 
with
\begin{equation*}
{\zeta_d}_m=\mathrm{mod}\left(\frac{{\xi_d}_m-a_d}{a_d-b_d},1\right)-\frac{1}{2},
\quad d=1,2,3
\end{equation*}
%where
%\begin{equation*}
%\mathrm{mod}(x,y)=x-y\left\lfloor\frac{x}{y}\right\rfloor
%\end{equation*}
and coefficients
\begin{equation*}
\hat f_{\boldsymbol k}=\hat\psi_{\boldsymbol
k}\prod_{d=1}^3\frac{\rme^{\pi\rmi(k_d-1-N_d/2)}}{\sqrt{b_d-a_d}}=\hat
\psi_{\boldsymbol k}
\frac{(-1)^{k_1+k_2+k_3-3-(N_1+N_2+N_3)/2}}{\prod_{d=1}^3\sqrt{b_d-a_d}}, 
\end{equation*}
where $\mathrm{mod}(x,y)$ is the usual remainder of the Euclidean division of
$x$ by $y$, 
$\mathrm{mod}(x,y)=x-y\left\lfloor x / y\right\rfloor$.

We first checked our MATLAB interface to NFFT by evaluating a truncated
series
of $64^3$ random complex Fourier coefficients at the regular grid 
$X_{\boldsymbol M}$ with $64^3$ points (for which the
inverse FFT is available) and compared the result in
infinity norm obtaining
\begin{verbatim}
NFFT

error_inf =

   5.3705e-14

Elapsed time is 2.267151 seconds.
\end{verbatim}
in \matlab{} and
\begin{verbatim}
NFFT
error_inf =    6.3161e-14
Elapsed time is 2.74056 seconds.
\end{verbatim}
in \octave{}.
Although the asymptotic cost of the NFFT is smaller than the evaluation
at the regular grid, for this number of coefficients and points of evaluation
NFFT turns out to be about 20 times slower than \verb+ndcovlt+.

In order to have an idea of the
computational cost in a real case usage, 
the evaluation of a series with $64^3$ coefficients
on $M=1000$ random points in $\Omega$ takes
\begin{verbatim}
one-loop
Elapsed time is 3.295020 seconds.
NFFT

error_inf =

   1.2296e-13

Elapsed time is 0.130843 seconds.
\end{verbatim}
in \matlab{} and
\begin{verbatim}
one-loop
Elapsed time is 3.95735 seconds.
NFFT
error_inf =    1.2488e-13
Elapsed time is 0.0957451 seconds.
\end{verbatim}
in \octave{}. Here we observe a speed-up of about 40 of the NFFT approach
over a straightforward implementation. The measured error is 
between the two evaluations. These tests can be found at the end of
the \verb+innft3.m+ file and can be run, in GNU Octave, by
\verb+demo infft3+.


%% References

\section{Application to quantum vortex reconnections}\label{sec:qvr}
Turbulence, ubiquitously present in nature, is dominated by reconnection of 
vortical structures.
Examples of reconnecting vortex tubes can be found in quantum 
turbulence~\cite{VIN2008,PL2011,BSS2014},
whose dynamics is properly described by the
Gross--Pitaevskii equation~(GPE)~\cite{P1961,G1963}
\begin{equation}\label{eq:GPE}
\frac{\partial \psi}{\partial t} = 
\Fra \rmi 2 \Lapl \psi + \Fra \rmi 2\left(1 - \abs{\psi}^2\right)\psi,
\end{equation}%
where $\psi$ is the complex wave function.
Quantum vortices are infinitesimally thin filaments of concentrated vorticity
in a unitary background density, 
$\rho(\boldsymbol x)=\abs{\psi(\boldsymbol x)}^2\to 1$ when 
$\abs{\boldsymbol x}\to \infty$.
On the vortex centerlines the density tends to zero and the phase of the wave 
function $\psi$ is not defined.
In the dimensionless units of equation~\eqref{eq:GPE}, the quantum of
circulation is $\Gamma=2\pi$ and the healing length, i.e. the lengthscale of the
core vortex over which reconnections occur, is $\xi=1$.
GPE conserves the (infinite) mass and the energy
\begin{equation}\label{eq:energy}
E=\frac{1}{2}\int \abs{\nabla \psi(\boldsymbol x)}^2\de \boldsymbol x+
\frac{1}{4}\int(1-\abs{\psi(\boldsymbol x)}^2)^2\de \boldsymbol x.
\end{equation}
Time splitting Fourier methods~\cite{KL93,ZCBB12,AZCPPB14,ZR15,CZ16x} are
normally used to compute the numerical solution of the GPE~\eqref{eq:GPE}.
%, given an initial condition where vortices are 
%imposed in the form of singular phase.
Because these methods rely on periodic boundary conditions
for the solutions restricted to a bounded physical domain,
initial conditions that are not periodic must be mirrored in the directions 
lacking periodicity~\cite{KL93}, with a consequent increase of the degrees of 
freedom and computational effort~\cite{CZ16x}.

Recent studies focusing on the topological details of quantum-vortex 
reconnections~\cite{ZR15} have emphasized the need for an accurate
description of the vortex centerline, which can be achieved by costly 
high-resolution numerical simulations of equation~\eqref{eq:GPE}.
One the other hand, it is possible to resort to more affordable approaches 
combined with an \emph{a posteriori} accurate evaluation of the solution on 
a finer grid, as proposed in Ref.~\cite{CZ16x}.

In order to exploit the second option, following~\cite{KL93}, we consider a 
fully three-dimensional reconnection originating from two perpendicular 
straight vortices, whose cross sections are two-dimensional vortices.
The wave function of a single two-dimensional vortex in the
$(s_1,s_2)$ plane and centered in $(0,0)$ is
$\rho(\sqrt{s_1^2+s_2^2})^{1/2}\rme^{\rmi \theta(s_1,s_2)}=
f(\sqrt{s_1^2+s_2^2})\rme^{\rmi \theta(s_1,s_2)}$, where
$f(\sqrt{s_1^2+s_2^2})=f(r)$ is a function to be determined whereas the phase is
$\theta(s_1,s_2)=\mathrm{atan2}(s_2,s_1)$.
By requiring the wave function to be the steady solution of equation~\eqref{eq:GPE}, 
we find~\cite{CZ16x} that $\rho(r)$ satisfies
\begin{equation}
\rho''+\frac{\rho'}{r}-\frac{(\rho')^2}{2\rho}
-\frac{2\rho}{r^2}+2(1-\rho)\rho=0,
\label{eq:GPErho}
\end{equation}
with boundary conditions $\rho(0) = 0$, $\rho(\infty)=1$. 
Instead of computing the numerical solution of this equation, it is possible to
resort to a high-order Pad\'e approximations of $\rho(r)$~\cite{CZ16x}.
It is known~\cite{B04,NW03} that diagonal 
Pad\'e approximations of $\rho(r)$ retain only
even degrees at both the numerator and denominator, that is
\begin{equation}
\rho(r) \approx
\rho_q(r)=
\frac{a_1r^2+a_2r^4+\cdots+a_qr^{2q}}{1+b_1 r^2+b_2r^4+\cdots+b_qr^{2q}}.
\end{equation}
The coefficients of a certain approximation $\rho_q(r)$ are computed by
substituting the analytic expressions $\rho_q(r)$, $\rho'_q(r)$ and 
$\rho''_q(r)$ in equation~\eqref{eq:GPErho} and by
nullifying the coefficients of the first $2q-1$ terms $r^{2k}$.
The choice $q=4$ leads to an algebraic equation of degree 8 for $a_1$ which 
can be solved numerically.
Once $a_1$ is known, all the other coefficients can be  computed analytically
(see~\cite{CZ16x} for the details).
Their expression is reported and used in the code file \verb+sf4pade.m+.
A straight vortex in a three-dimensional domain can be obtained by
the extrusion of the above two-dimensional wave function along the vortex
center line.
A nontrivial initial condition generated by the superimposition of 
multiple straight vortices is simply the product of their 
wave functions.
\subsection{Numerical discretization}
After restricting the unbounded domain $\RR^3$ 
to the computational
domain $\Omega=\prod_{d=1}^3[a_d,b_d)$ in which the initial
solution is periodic, equation~\eqref{eq:GPE} can be split into
the \emph{kinetic} and \emph{potential} parts
\begin{subequations}
\begin{align}
\frac{\partial u}{\partial t}&=\frac{\rmi}{2}\Lapl u\label{eq:GPEkin}\\
\frac{\partial v}{\partial
t}&=\frac{\rmi}{2}\left(1-\abs{v}^2\right)v,\label{eq:GPEpot}
\end{align}
\end{subequations}
and the Time Splitting 
pseudoSPectral~(TSSP) approach can be employed, as done in~\cite{CZ16x}.
Equation~\eqref{eq:GPEkin} is solved exactly in time within the Fourier 
spectral space,
%While smooth solutions fastly decaying to zero ensure the spectral order of 
%convergence, when simulating the dynamics of vortex solutions in a non-zero
%background (here $\lim_{\abs{x}\to\infty}\abs{\psi(t,x)}=1$)
%some issues such as the low regularity of the solution and the lack of 
%periodicity at the boundaries must be taken care of (see~\cite{CZ16x} for a
%thorough discussion).
whereas equation~\eqref{eq:GPEpot} is solved exactly owing to the fact that 
$\abs{v}$ is preserved by the equation. Therefore,
\begin{equation}
v(\tau,\boldsymbol x)=\exp\left(\frac{\tau\rmi}{2}\left(1-\abs{v(0,\boldsymbol x)}^2\right)\right)v(0,\boldsymbol x)
\end{equation}
for any $\boldsymbol x$ in the spatial domain. 
By introducing 
$\rme^{\tau\mathcal{A}}u_n(\boldsymbol x)$ and 
$\rme^{\tau\mathcal{B}(v_n(\boldsymbol x))}v_n(\boldsymbol x)$ to denote the two 
partial numerical 
solutions, the numerical approximation $\psi_{n+1}(\boldsymbol x)$ of $\psi(t_{n+1},\boldsymbol x)$ at
time $t_{n+1}=(n+1)\tau$ is recovered by the so-called Strang splitting
\begin{equation*}
\begin{aligned}
\psi_{n+1/2}(\boldsymbol x)&=\rme^{\tau\mathcal{A}}
\rme^{\frac{\tau}{2}\mathcal{B}(\psi_n(\boldsymbol x))}\psi_n(\boldsymbol x)\\
\psi_{n+1}(\boldsymbol x)&=\rme^{\frac{\tau}{2}\mathcal{B}(\psi_{n+1/2}(\boldsymbol x))}\psi_{n+1/2}(\boldsymbol x).
\end{aligned}
\end{equation*}
Strang splitting preserves the discrete finite mass in the computational
domain $\Omega$ and is second order accurate in time. We refer 
the reader to~\cite{TCN09} for higher-order time splitting methods.

\subsection{Other applications of the NFFT tool}
The Gross--Pitaevskii equation is a model not only for superfluids but
also for Bose--Einstein condensates (see~\cite{BC13} for a review). 
In the second framework, the typical formulation is
\begin{equation*}
\rmi \frac{\partial \psi}{\partial t}=-\frac{1}{2}\nabla^2\psi+
V\psi+\beta \lvert\psi\rvert^{2\sigma}\psi,
\end{equation*}
where $V\colon \RR^3\to\RR$ is a scalar potential, $\beta$ a real constant
and $\sigma >0$. 
In this case the corresponding energy is
\begin{equation*}
E=\frac{1}{2}\int \lvert\nabla \psi(\boldsymbol x)\rvert^2\rmd \boldsymbol x+\int V\lvert\psi(\boldsymbol x)\rvert^2\rmd \boldsymbol x
+\frac{\beta}{\sigma+1}\int \lvert\psi(\boldsymbol x)\rvert^{2\sigma+2}\rmd
\boldsymbol x
\end{equation*}
and the Strang splitting method described above can still be applied without 
any modification.
We notice that space discretizations which are not regular
(see, for instance, \cite{CZ16x} for nonuniform finite differences and
\cite{TA12} for finite elements) provide results that are difficult to compare
with those obtained via pseudospectral approaches, which are available only on
regular grids.
INFFTM allows the evaluation at arbitrary rectilinear grids and sets of 
arbitrary points making the comparison of these results possible.

Another interesting application where NFFT is a valuable tool 
is the so called \emph{magnetic} Schr\"odinger equation
\begin{equation*}
\rmi\frac{\partial \psi}{\partial t}=\frac{1}{2}(\rmi\nabla+A)^2\psi+V\psi,
\end{equation*}
where $A\colon \RR^3\to\RR^3$ is the vector potential which can be chosen
divergence free owing to Coulomb's gauge.
Besides the kinetic and the potential parts, 
the \emph{advection} part
\begin{equation*}
\frac{\partial w}{\partial t}=A\cdot \nabla w
\end{equation*}
has to be considered 
and then combined with the others in a splitting scheme.
The advection part can be solved,
for instance, by the characteristics method and the value of 
$w(\tau,\boldsymbol x)$ at the departure point of the characteristics
can be recovered by NFFT. We refer to~\cite{COP15x} for further details.
\section{Description of the programs}\label{sec:prog}
%We describe here the main functions used by the code. 
On developing the code, we realized that some functions naturally apply
to any space dimension. On the contrary, others are specific for the 
three-dimensional case, which is the object of the present work. Therefore,
we used the following
convention: function names ending in `3' are specific and 
for the three-dimensional case only, whereas the others can work in any
space dimension. 
The only exception is \verb+igridftn+ which calls \verb+ndcovlt+, originally 
developed by Jaroslav Hajek and not designed for the trivial 
one-dimensional case (see \S~\ref{sec:fse} and \ref{sec:aux}).
In what follows we describe only the implementation in three dimensions.

As written in the \verb+README+ file, before using the package, NFFT has
to be installed. We refer to \ref{sec:NFFTLinux} for the
instructions on the installation in a Linux environment. After that,
%the package can be decompressed and
%the variable \verb+havenfft+ has to be set to \verb+true+ in the file
%\verb+nfftpath.m+. If the GNU Octave \verb+nfftpkg-0.0.3.tar.gz+ package
%was installed, the variable \verb+havenfftpkg+ has to be set to \verb+true+,
%too. Otherwise, 
the correct path to the NFFT library has to be given
in the file \verb+nfftpath.m+. If the NFFT library is not installed,
the package will work anyway, but the evaluation of a three-dimensional 
truncated
Fourier
series at a set of arbitrary points will be extremely slow.
\subsection{Functions for Fourier series evaluation}\label{sec:fse}
%As written in the \verb+README+ file, before using the package, NFFT has
%to be installed. We refer to section~\ref{sec:NFFTLinux} for the
%instructions on the installation on a Linux environment. After that,
%the package can be decompressed and
%the variable \verb+havenfft+ has to be set to \verb+true+ in the file
%\verb+nfftpath.m+. If the GNU Octave \verb+nfftpkg-0.0.3.tar.gz+ package
%was installed, the variable \verb+havenfftpkg+ has to be set to \verb+true+,
%too. Otherwise, the correct path to the NFFT library has to be given
%in the file \verb+nfftpath.m+. If the NFFT library is not installed,
%the package will work anyway, but the evaluation of a 3d Fourier
%series at a set of arbitrary points will be extremely slow.
%
The two main functions are \verb+igridftn+ and \verb+infft3+. They implement
the evaluation of the truncated Fourier series~\eqref{eq:fs} at a rectilinear
grid (\verb+ndgrid+ format) and at an arbitrary set of points, respectively. 
The calls are similar
\begin{verbatim}
psi = igridftn(psihat,a,b,y)
psi = infft3(psihat,a,b,Xi)
\end{verbatim}
\verb+psihat+ being the three-dimensional array of Fourier coefficients,
\verb+a+ and \verb+b+ the limits of the physical domain $\Omega$ (in the form
\verb+[a(1),a(2),a(3)]+ and 
\verb+[b(1),b(2),b(3)]+), \verb+y+ a cell array containing
in the column vector \verb+y{d}+ the $d$-th projection of the points and 
\verb+Xi+ a 
two-dimensional array containing in the $d$-th row the $d$-th component of
the points.

The simple function 
\begin{verbatim}
plotiso3(x,data,iso)
\end{verbatim}
invokes the MATLAB program \verb+isosurface+ to plot the isosurface of level
\verb+iso+ of the real input \verb+data+ corresponding to \verb+ndgrid{x{1:3}}+.
Since \verb+isosurface+ in \matlab{} requires the data in \verb+meshgrid+ 
format, \verb+plotiso3+ performs the permutation
\begin{verbatim}
data = permute(data,[2,1,3]); 
\end{verbatim}
\subsection{Functions for superfluid simulation by GPE}
As described in~\S~\ref{sec:qvr}, here we focus on the particular application to
quantum vortex reconnections.
The time integration of GPE is carried out by the main function \verb+sfrun+.
Given an initial solution as a function of \verb+x{1:3}+, \verb+sfrun+
first computes some preliminary quantities (\verb+sfpregpe+), such as 
$\Lambda_{\boldsymbol k,d}$ (see eq.~\eqref{eq:Lambda}),
then computes initial and final mass and energy of the system 
(\verb+sfEm+, see eq.~\eqref{eq:energy}), 
and finally it performs time integration
by Strang splitting method (\verb+sfgpe+) and store the structure
\verb+sf+ of the solution in a MATLAB \verb+'-v6'+ format file
at \verb|nsteps+1| equally distributed time steps. 
The structure \verb+sf+
contains the fields \verb+pdb+ (a row vector of length six consisting of
the physical domain boundaries), \verb+psipdb+ (a complex 3d-array of the 
values of the
wave function at the grid in the physical domain), \verb+mirror+ (a row
vector of length three for the mirroring flags) and \verb+t+ (the
simulation time). 
From the structure \verb+sf+ it is possible to
recover the Fourier coefficients of the solution \verb+psipdb+ by invoking
the function
\begin{verbatim}
[psihat,a,b] = sf2psihat(sf)
\end{verbatim}
All the previous functions work in any space dimension. The functions
\verb+sfsvl3+ and \verb+sfic3+ generate respectively a single
straight vortex in a three-dimensional domain and the superimposition
of multiple vortices.
The function \verb+sfview3+ simply
extracts the grid points and the density of the wave function
from the structure \verb+sf+ and plots a given isosurface level through the
function \verb+plotiso3+.
\subsection{Evaluation within a vortex tube}\label{sec:sftubecoll}
The study of vortex reconnections in quantum fluids requires
high spatial resolution in order to extract the vortex centerlines with enough
accuracy, and this is especially true in the neighborhood of the 
reconnection event (see~\S~\ref{sec:qvr}).
Instead of evaluating the physical solution at a finer rectilinear grid 
within the 
whole physical domain, it is more convenient to evaluate the solution only
within vortex tubes, i.e. where high resolution is really needed.

Function \verb+sftubeeval3+ has been designed especially for this purpose.
Its input arguments are the structure \verb+sf+, which defines completely 
$\psi$ on an equispaced grid in the physical domain, and \verb+rhobar+,
a vector containing the values of the density $\rho$ that define the
vortex tubes.
For example, if \verb+rhobar=0.2+ (a single value), then function
\verb+sftubeeval3+ first extracts points from $X_{\boldsymbol N}$ 
for which $\rho \le 0.2$.
Then, if $\boldsymbol \xi_{m}$ denotes the $m$-th point within the 
vortex tube (corresponding to a certain $\boldsymbol x_{\boldsymbol n}\in X_{\boldsymbol N}$)
and $h_d$ the step-size of $X_{\boldsymbol N}$
in direction~$d$, a small 
regular grid of step-size $h_d/3$ made of only 27 points 
centered in $\boldsymbol \xi_{m}$ is generated for each~$m$.
Finally, \verb+sftubeeval3+  returns the new set of points and 
$\rho = |\psi|^2$ evaluated at these points by NFFT.
In order to retrieve smaller vortex tubes containing enough points, the input
\verb+rhobar+ should be a vector.
In this case the process described above is repeated up to the last value
of $\rho$ and the output of \verb+sftubeeval3+  is the set of points $\Xi_M$
on successive refined grids
for which $\rho \le$~\verb+rhobar(end)+, together with their corresponding 
values of $\rho$.
\subsection{Drivers}
The two drivers \verb+sfdrv3+ and \verb+evaldrv3+
were written for the convenience of the user, as they perform the numerical 
simulation and the visualizations exactly as described in the next section.
\section{Numerical experiments}\label{sec:numexp}
\begin{figure}[!ht]
\centering
\includegraphics[scale=0.6]{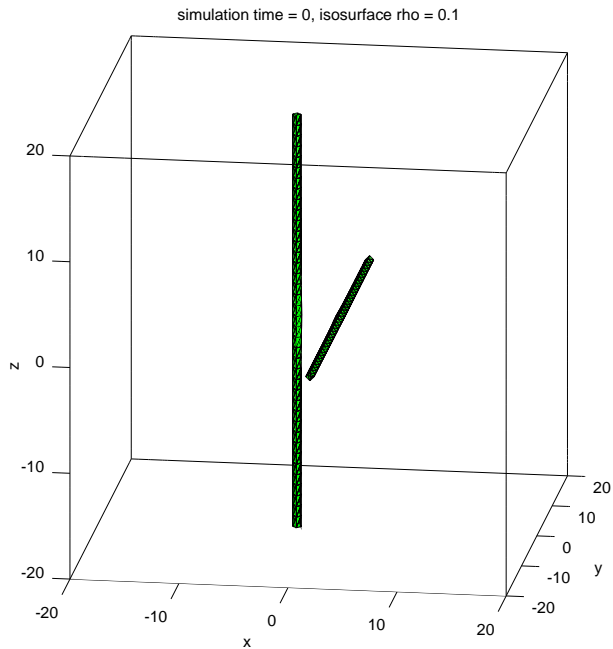}
\caption{Isosurface level 0.1 for the density of the initial solution at $X_{\boldsymbol N_1}$.}
\label{fig:inisol}
\end{figure}
We solve the GPE equation~\eqref{eq:GPE} in the physical domain 
$[-20,20]^3$. The initial solution is given by the superimposition of two
straight vortices, passing through the points $(2,0,0)$ and $(-2,0,0)$
and oriented as $(0,1,0)$ and $(0,0,1)$, respectively. 
In order to make this initial condition periodic %at least in the values 
at the boundaries,
the computational domain has to be set to $\Omega=[-20,60)^3$ and the
initial solution has to be mirrored along the three directions.
This can be accomplished by setting the field \verb+mirror+ to
\verb+[true,true,true]+ in the structure \verb+sf+ associated to the initial
solution.
While ensuring the periodicity of the solution, mirroring does not force the
periodicity of the derivatives. 
In the computational domain we select $\boldsymbol N_1=(80,80,80)$, yielding
a (coarse) regular grid $X_{\boldsymbol N_1}$ with a constant space step size 
of 1 along each direction.
The solution is computed up to the final time $T=20$ with 200 time steps.
The initial solution, in the original physical domain, is shown in 
Figure~\ref{fig:inisol} by plotting the isosurface $\rho=0.1$ extracted from 
the original data at the regular grid $X_{\boldsymbol N_1}$,
through the function \verb+sfview3+.
\begin{figure}[!ht]
\centering
\includegraphics[scale=0.45]{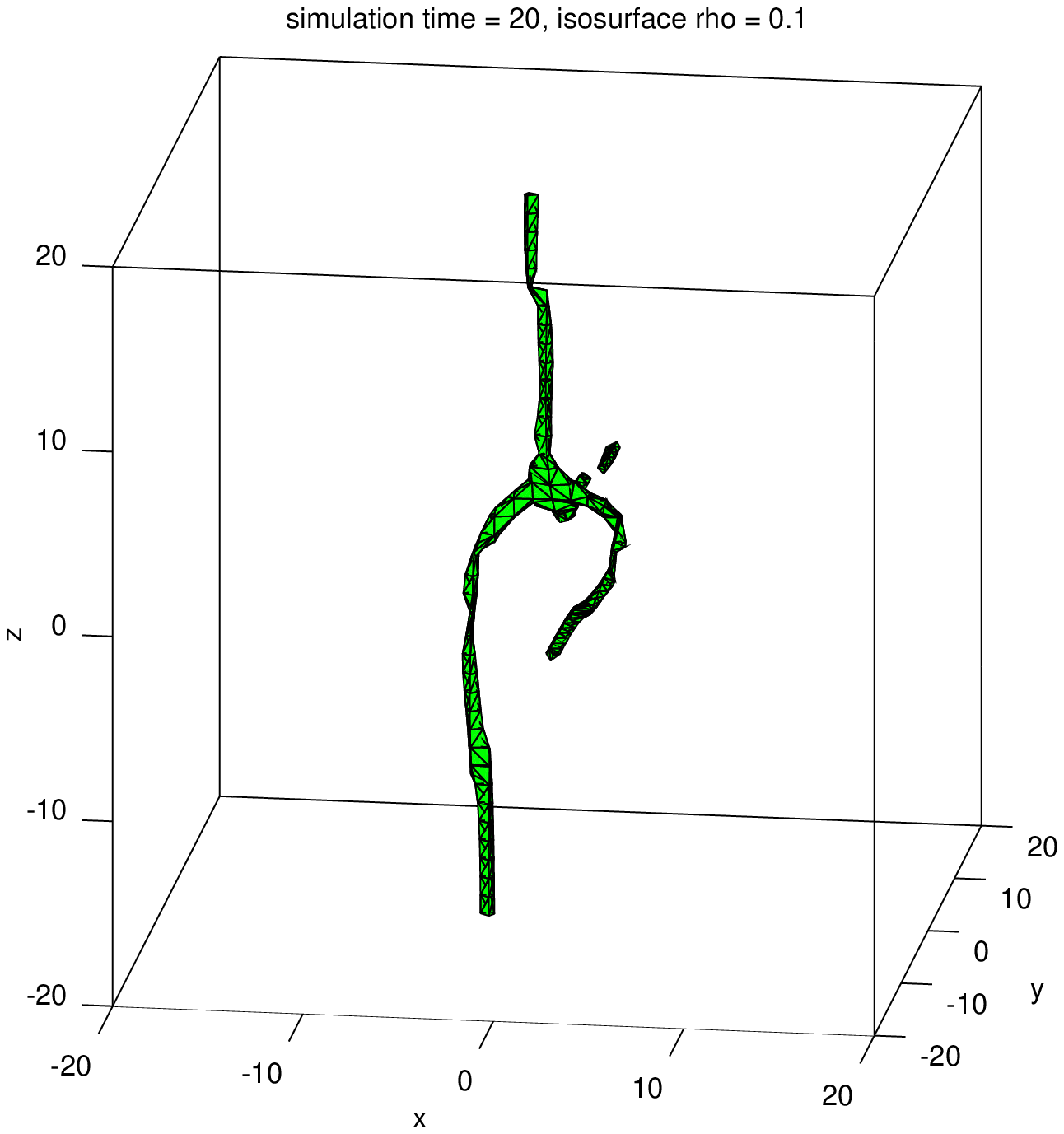}
\includegraphics[scale=0.45]{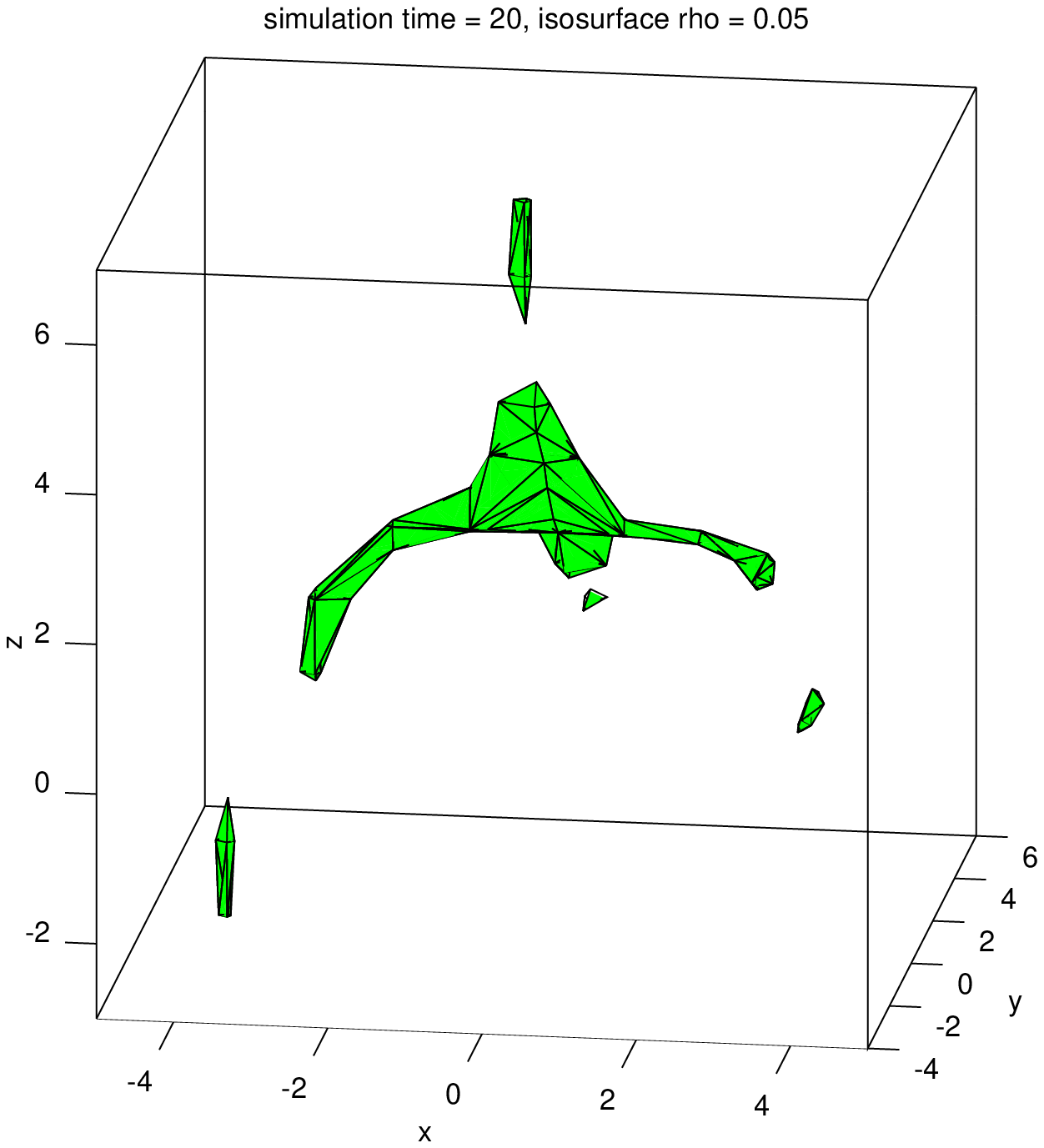}
\caption{Isosurface level 0.1 (left) 
and zoom of the isosurface level 0.05 (right) for the density of the final solution at $X_{\boldsymbol N_1}$.}
\label{fig:finsol}
\end{figure}

\begin{figure}[!ht]
\centering
\includegraphics[scale=0.45]{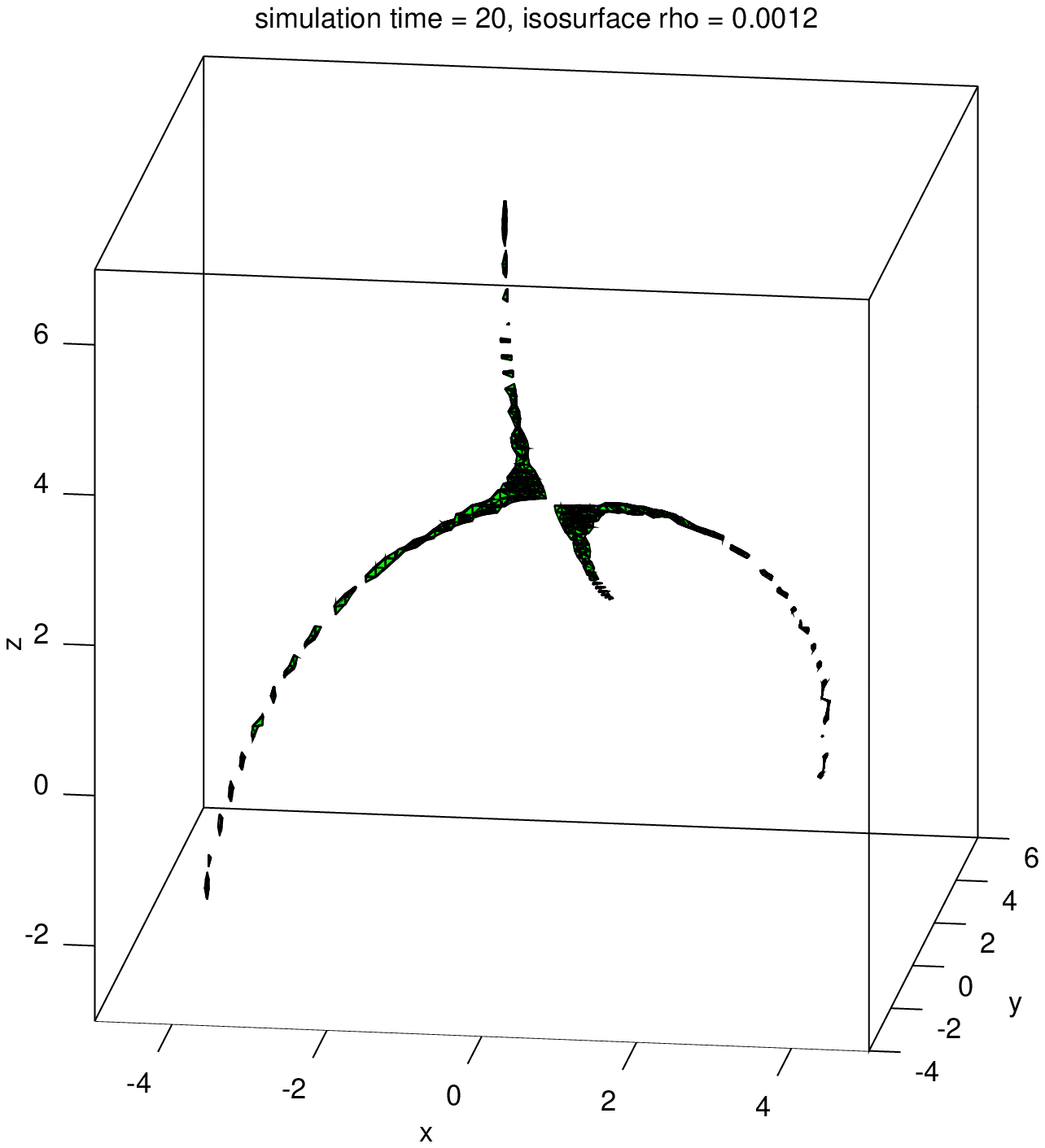}
\includegraphics[scale=0.45]{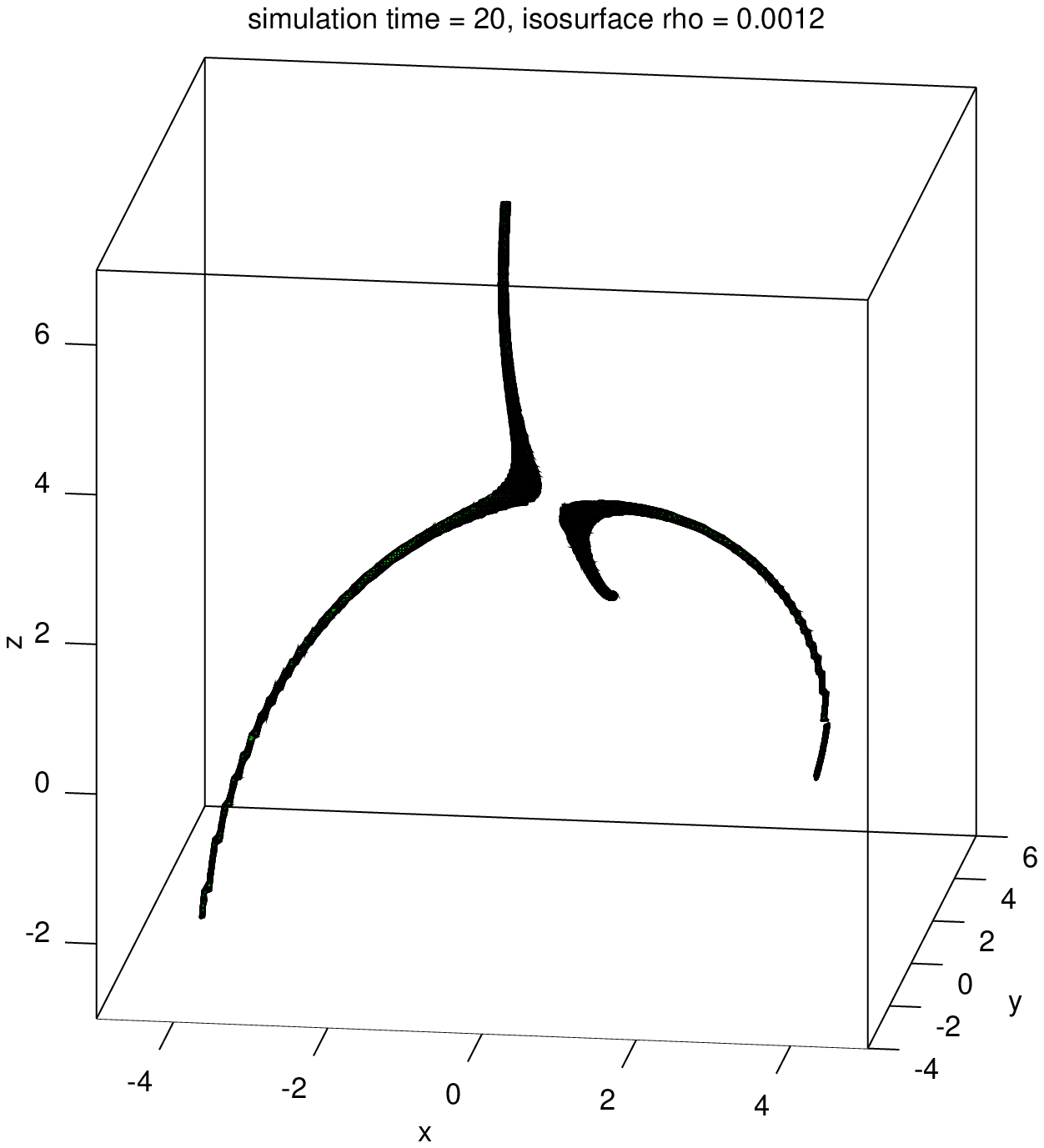}
\caption{Isosurface level 0.0012 for the density of the
final solution evaluated at $Y_{\boldsymbol M_1}$
(left, equispaced grid) and at $Y_{\boldsymbol M_2}$ (right, nonequispaced grid).}
\label{fig:finsolcoll}
\end{figure}
The solution at the final time $T=20$ is reported in Figure~\ref{fig:finsol}.
The left plot shows the isosurface $\rho=0.1$ in the whole physical domain, 
whereas the right plot shows a zoom at the isolevel $\rho=0.05$, which should 
guarantee a better description of the vortex centerlines ($\rho\to 0$ therein).
Clearly, none of the plots in Figure~\ref{fig:finsol} allows to discriminate
whether the reconnection has occurred or not.
Moreover, reducing the isolevel of $\rho$ makes things worse in that vortex
tubes appear disconnected due to the low spatial resolution 
characterizing the original data. 

In order to increase the details, we evaluate the solution at a finer
Cartesian equispaced grid $Y_{\boldsymbol M_1}$, with 
$\boldsymbol M_1=(321,321,321)$ in the physical domain $[-20,20]^3$ by
the function \verb+igridftn+.
By extracting the isosurface corresponding to $\rho=0.0012$, 
the vortex tubes become much better defined (see Figure~\ref{fig:finsolcoll}, 
left) clearly indicating that a reconnection has occurred. 
Selecting the same isolevel for the original data at
$X_{\boldsymbol N_1}$ yields an almost empty plot.
Since our interest is in the neighborhood of the reconnection, instead of 
evaluating the solution at equally-spaced points, it is more convenient
to evaluate the solution at a coarser \emph{nonequispaced} rectilinear grid 
$Y_{\boldsymbol M_2}$, with $\boldsymbol M_2=(281,281,281)$ points 
denser around the origin, always by \verb+igridftn+.
The isosurface $\rho=0.0012$ (see Figure~\ref{fig:finsolcoll}, right)
provides a much better result than the equispaced case 
(Figure~\ref{fig:finsolcoll}, left)
in terms of clear vortex cores, which now appear completely connected.
\begin{figure}[!ht]
\centering
\includegraphics[scale=0.45]{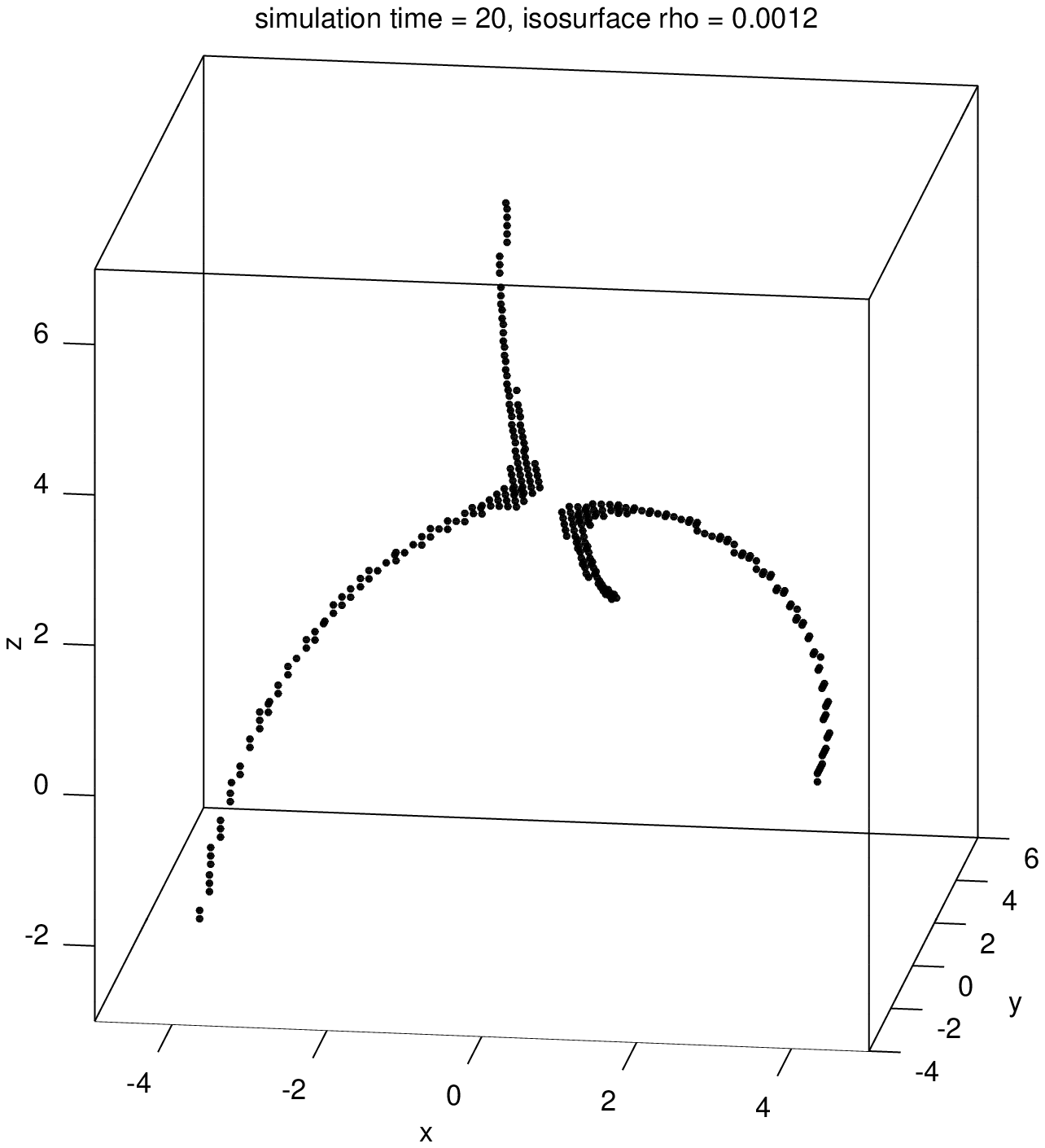}
\includegraphics[scale=0.45]{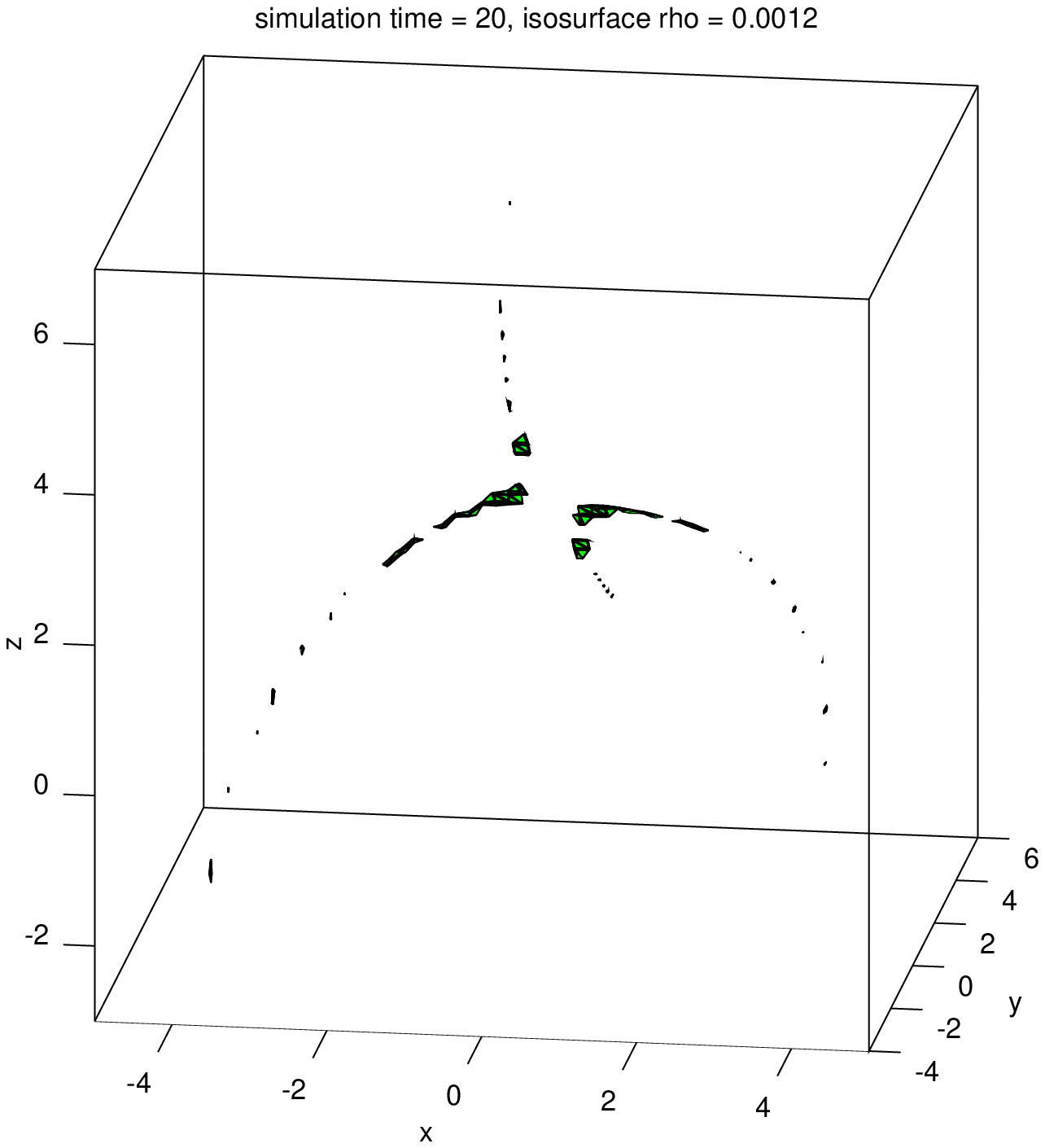}
\caption{Isosurface level 0.0012 for the density of the
final solution evaluated at $\Xi_M$ (left) and
for the density of the final solution computed and evaluated
at $X_{\boldsymbol N_2}$, with $\boldsymbol N_2=(228,228,228)$ (right).}
\label{fig:finsolnfft}
\end{figure}

Finally, by employing the strategy described in~\S~\ref{sec:sftubecoll} with the
sequence \verb+rhobar=[0.2,0.05]+ and then plotting points corresponding to
 $\rho\le 0.0012$, we obtain the left plot in 
Figure~\ref{fig:finsolnfft}. The set $\Xi_M$ obtained by \verb+sftubeeval3+
has $M=32022$ points and the number of points corresponding
to $\rho\le 0.0012$ is 809.
%is possible to observe the reconnection by a evaluation through NFFT
%at a set $\Xi_M$ of $M=32022$ points, iso 0.0012 corresponds to 809 
Similar vortex tubes can be obtained without evaluation at finer grids only by
resorting to high-resolution simulations.
An example is reported in
the right plot of Figure~\ref{fig:finsolnfft}, which shows the isosurface $\rho=0.0012$
for the solution at the original regular grid $X_{{\boldsymbol N}_2}$
with $\boldsymbol N_2=(228,228,228)$ and 1000 time steps.

The script to run the GPE simulation with $\boldsymbol N_1$
is \verb+sfdrv3+ (CPU time about 14~s),
whereas the script to perform evaluation is
\verb+evaldrv3+ (CPU time about 17~s with NFFT installed, about 247~s 
without).
The simulation with
$\boldsymbol N_2$ took about 420 minutes.
\section{Conclusions}
We have developed the package INFFTM for the fast evaluation
of three-dimensional truncated Fourier series at general rectilinear
grids and sets of arbitrary points. The two main functions,
\verb+igridftn+ and \verb+infft3+, are written in plain MATLAB language, work
in Matlab$^\text{\textregistered}$ and GNU Octave and are based on two
efficient, although not widespread, tools, namely \verb+ndcovlt+
by J.~Hajek and NFFT by J.~Keiner, S.~Kunis, and D.~Potts. 
We have demonstrated the effectiveness of \verb+igridftn+ and \verb+infft3+ 
in the framework of quantum vortex reconnections. 
A proper post-processing of the numerical data
obtained by running a cheap simulation of the vortex dynamics modeled 
by the Gross--Pitaevskii equation provides details on the reconnecting vortices
that are comparable to costly high-resolution simulations.
These promising results highlight the potential of INFFTM to become a 
standard MATLAB library for applications involving Fourier series approximation.
\appendix
\section{Installation of NFFT in a Linux environment}\label{sec:NFFTLinux}
The necessary information for installation of NFFT is available in
the file \verb+README+. Here we briefly summarize the procedure.

On writing this paper, the latest release of NFFT was
\verb+nfft-3.3.1.tar.gz+\footnote{Available at \texttt{https://www-user.tu-chemnitz.de/$\sim$potts/nfft/}.}.
It can be built in the usual way (\verb+./configure+ and \verb+make+) and,
in order to compile the Matlab$^\text{\textregistered}$ 
mex interface, it must to be  configured
%\footnote{With Matlab$^\text{\textregistered}$ R2014b on a Linux/64bit system we had to add the option \texttt{LDFLAGS=-Wl,-rpath,/path/to/matlab/bin/glnxa64} to \texttt{./configure}.} 
by
\begin{verbatim}
./configure --with-matlab=/path/to/matlab --enable-openmp
\end{verbatim}
Then, the correct path to \verb+nfft-3.3.1/matlab/nfft+ has to be given in
the \verb+nfftpath.m+ file. The path to Matlab can be obtained by the command
\begin{verbatim}
matlab -n
\end{verbatim}

The installation in GNU Octave is not difficult but requires some
modifications to the \verb+configure+ script.
For the user convenience,
in the \verb+aux+ folder we provide the file
\verb+nfft-3.3.1-octave.patch+
which has to be put in the \verb+nfft-3.3.1+ folder and then
applied by typing the command
\begin{verbatim}
patch -p1 < nfft-3.3.1-octave.patch
\end{verbatim}
from within the folder \verb+nfft-3.3.1+. Finally, the configuration is made
by
\begin{verbatim}
./configure --with-octave=/path/to/octave/headers --enable-openmp
\end{verbatim}
and then \verb+make+.
\verb+/path/to/octave/headers+ is the folder containing
\verb+octave.h+ and can be recovered by giving the shell command
\begin{verbatim}
mkoctfile -p OCTINCLUDEDIR
\end{verbatim}

Instead of manually patching the original sources of \verb+nfft-3.3.1+,
it is possible to install the
package  \verb+nfftpkg-0.0.4.tar.gz+
that we built for GNU Octave.
It is contained in the \verb+aux+ folder\footnote{The patches and
the packages for GNU Octave are also available at corresponding author's homepage
\texttt{http://profs.scienze.univr.it/caliari/software.htm}.}
and can be installed in the usual way under GNU Octave
\begin{verbatim}
octave:1> pkg install nfftpkg-0.0.4
\end{verbatim}
%and loaded when needed
%\begin{verbatim}
%octave:2> pkg load nfftpkg
%\end{verbatim}
%If NFFT is installed in this way, there is no need to set the
%path in the \verb+nfftpath.m+ file 
%and \verb+havenfftpkg+ should be set to true therein. 
In order to 
check the correct installation of the NFFT library it is possible to
run the test at the end of the \verb+nfftpath.m+ file or, in GNU Octave,
to run \verb+demo nfftpath+. The demo (provided by the original NFFT library) 
requires about 20 s and has to be considered passed if the string
\verb+A two dimensional example+ appears. The presence of
some \verb+NaN+ values in the output is \emph{not} a symptom of a failure.
\section{Auxiliary files and workarounds}\label{sec:aux}
%Invalid MEX-file '/usr/local/src/nfft-3.3.1-matlab/matlab/nfft/nfftmex.mexa64': %dlopen:
%cannot load any more object with static TLS
Versions of GNU Octave before 4.0.0 have no \verb+flip+ function, which
is required by the code and distributed in the \verb+aux+ folder. % of the code. 
GNU Octave 4.0.0. has a 
bug\footnote{Bug no.~\#45207.}
with \verb+fftshift+ and \verb+ifftshift+ not working on three-dimensional
arrays. Patched working functions are distributed in the \verb+aux+ folder.

\matlab{} has a bug\footnote{Bug no.~961694.} preventing, form time to time,
to load \verb+nfftmex.mexa64+. The workaround is to load the library as soon
as \matlab{} is started. This can be achieved, for instance, by running the script \verb+aux/mlloadnfft+.

In the folder \verb+aux+ we provide also the functions \verb+igridft2+, 
\verb+infft+, and \verb+infft2+. Although \verb+igridftn+ and 
\verb+ndconvlt+ can work in two dimensions, in order to evaluate
a two-dimensional truncated 
Fourier series at a rectilinear grid it is much simpler to use
a double matrix-matrix product as done in \verb+igridft2.m+.
The functions \verb+infft+ and \verb+infft2+
apply NFFT to a one-dimensional and to a two-dimensional array
of Fourier coefficients, respectively. A demonstration of their usage
can be found at the end of the files and can be run, in GNU Octave, by
\verb+demo infft+ and \verb+demo infft2+.
We notice that there is no need for a specific function \verb+igridft1+ because
the evaluation of a one-dimensional truncated Fourier series at an
arbitrary set $\Xi_M$ of points is straightforwardly obtained by 
the matrix-vector product
\begin{verbatim}
E = exp(2*pi*1i * (Xi(:) - a) * (-N/2:N/2 - 1) / (b - a)) / ...
        sqrt(b - a);
psi = E * psihat;
\end{verbatim}
whose computational cost is $\mathcal{O}(NM)$. However, 
the one-dimensional \verb+infft+ might be more convenient for large $M$
as its computational cost is $\mathcal{O}(N\log N+M\abs{\log \varepsilon})$
and, in general, it is also more accurate.






% \end{thebibliography}

\end{document}